\documentclass{amsart}
\usepackage{amssymb} 
\usepackage{amsthm, multirow, subfigure}
\usepackage[all]{xy}
\usepackage{graphicx}

\newcommand{\Z}{\mathbb{Z}}
\newtheorem{theorem}{Theorem}[section]
\newtheorem*{theorem*}{Theorem}
\newtheorem{lemma}[theorem]{Lemma}
\newtheorem{proposition}[theorem]{Proposition}
\newtheorem{corollary}[theorem]{Corollary}
\theoremstyle{definition}
\newtheorem{example}[theorem]{Example}
\theoremstyle{remark}
\newtheorem*{remark}{Remark}
\title{Substitutions and $\frac{1}{2}$-Discrepancy of $\{n \theta + x\}$}
\author{David Ralston}
\email{ralston.david.s@gmail.com}
\address{Ben Gurion University, Department of Mathematics \\ POB 653 \\ Beer Sheva 84105 \\ ISRAEL}
\date{\today}
\begin{document}
\begin{abstract}
The sequence of $1/2$-discrepancy sums of $\{x + i \theta \bmod 1\}$ is realized through a sequence of substitutions on an alphabet of three symbols; particular attention is paid to $x=0$.  The first application is to show that any asymptotic growth rate of the discrepancy sums not trivially forbidden may be achieved.  A second application is to show that for badly approximable $\theta$ and any $x$ the range of values taken over $i=0,1,\ldots n-1$ is asymptotically similar to $\log(n)$, a stronger conclusion than given by the Denjoy-Koksma inequality.
\end{abstract}
\maketitle
\section{Introduction}

Given an irrational $\theta$ and some $x \in [0,1)=S^1$ (all addition in $S^1$ is taken modulo one), let
\begin{equation}\label{eqn - function f}
f(x) = \chi_{[0,1/2)}(x) - \chi_{[1/2,1)}(x).\end{equation}
With $\theta$ fixed, the \textit{$1/2$-discrepancy sums} of the sequence $\{x+i\theta\}$ are given by
\[S_n(x) = \sum_{i=0}^{n-1}f(x+i\theta).\]
Two results are classical in this setting, for any irrational $\theta$ and for all $x$:
\begin{equation}\label{eqn - unbounded but not linear}
S_n(x) \in o(n), \quad S_n(x) \notin O(1).
\end{equation}
The first restriction is due to unique ergodicity of the underlying rotation, and the second is a theorem of Kesten \cite{MR0209253}.

We will use standard continued fraction notation; partial quotients are denoted $a_i(\theta)$, and convergents are denoted $p_i(\theta)/q_i(\theta)$.  When $\theta$ is clear from context we will simply write $a_i$, $p_i$ and $q_i$.  The distance from $x$ to the nearest integer is denoted $\| x \|$. As $\theta \in (0,1)$ without loss of generality, we will assume that $a_0(\theta)=0$ and omit this term, writing simply
\[\theta = [a_1,a_2,a_3,\ldots]=\cfrac{1}{a_1+\cfrac{1}{a_2+\cfrac{1}{a_3+\ddots}}}.\]  All necessary background in continued fractions may be found in \cite{MR1451873}.  The \textit{Gauss map} will be denoted by $\gamma$, and acts as the non-invertible shift on the sequence of partial quotients:
\begin{equation}\label{eqn - gauss map}
\gamma(\theta) = \frac{1}{\theta} \bmod 1, \quad \gamma([a_1,a_2,\ldots])=[a_2,a_3,\ldots].
\end{equation}

Our goal is to investigate what behavior is possible for the sequence $S_n(x)$ within the constraints of \eqref{eqn - unbounded but not linear}.  Because the sequence $S_n$ is not monotone, however, it will be more convenient to consider the following sequences, which track the \textit{maximal} and \textit{minimal} discrepancies, as well as the \textit{range} of values taken:
\begin{gather}
M_n(x) = \max\{S_i(x): i=1,\ldots,n-1\},\\
m_n(x) = \min\{S_i(x): i=1,\ldots,n-1\},\\
\rho_n(x) = M_n(x) - m_n(x)+1.
\end{gather}
It is worth clarifying that $m_n$ is taken as a minimum over \textit{integers}, and as such can generally be expected to be negative.  It is a matter of later convenience that $i=0$ is not considered: for example, $M_1(0)=m_1(0)=S_1(0)=1$.

We will develop a renormalization procedure through which the sequence of values $f(x+i\theta)$ can be determined from a sequence of substitutions.  Let $\theta<1/2$ and $A=[0,1/2)$, $B=[1/2,1-\theta)$, $C=[1-\theta,1)$.  If we wish to change which interval certain endpoints belong to (for example, if we wish for $A$ to be closed and $B$ to be open), we will say that we make a \textit{change of endpoints} of the intervals $A$, $B$, and $C$.  Our central result is the following:
\begin{theorem}\label{theorem - substitution sequence}
Given any irrational $\theta$ and any $x \in [0,1)$, there is a sequence of words $\omega_i$ (some of which may be empty) and substitutions $\sigma_i$ (infinitely many are not identity) both defined on the alphabet $\{A,B,C\}$, given by a dynamic process depending on $x$ and $\theta$, such that the infinite word given by
\begin{equation}\label{eqn - arbitrary substitution point}\omega_0 \sigma_0\left( \omega_1 \sigma_1 \left( \omega_2 \sigma_2 ( \ldots) \right) \right)\end{equation} encodes the orbit of $x$ up to at most two errors.  Alternately, the coding is exact up to a change of endpoints of the intervals $A$, $B$ and $C$.  The dependence of $\sigma_i$ on $\theta$ and $\omega_i$ on $(x,\theta)$ is explicit.

There is one special point $x(\theta)$ for which \textit{all} $\omega_i$ may be taken to be the empty word, in which case the infinite word
\begin{equation}\label{eqn - special substitution point}\lim_{n \rightarrow \infty}\left(\sigma_0 \circ \sigma_1 \circ \cdots \circ \sigma_{n-1}\right)(\omega)\end{equation} will encode the orbit of $x(\theta)$ regardless of the choice of nonempty word $\omega$.  The orbit of zero can alternately be determined by
\begin{equation}\label{eqn - zero substitutions}\lim_{n \rightarrow \infty}\left( \sigma'_0 \circ \sigma'_1 \circ \cdots \circ \sigma'_{n-1}\right)(\omega'_{n-1}),\end{equation} where $\sigma'_n$ are either substitutions or a different map.  This distinction and the word $\omega'_n$ are explicitly presented.\end{theorem}

We will include some remarks regarding the point $x(\theta)$ (including a complete characterization of those $\theta$ for which $x(\theta)=0$ in Proposition \ref{proposition - heavy theta have zero}), as well as proving that the sequence of substitutions $\sigma_i$ is eventually periodic if and only if $\theta$ is a quadratic surd (Proposition \ref{proposition - quadratics are periodic}).

As $(0,1/2)\subset A$ and $(1/2,1)\subset(B \cup C)$, any change of endpoints is completely irrelevant to the asymptotic growth rates of $M_n(x)$, $m_n(x)$, and $\rho_n(x)$.  While Theorem \ref{theorem - substitution sequence} provides a way to produce the orbit of an arbitrary point, computation of the words $\omega_i$ is a nontrivial task.  However, for the special point $x(\theta)$ and for $0$, the process is much simpler.  We will show that given any growth condition that does not violate \eqref{eqn - unbounded but not linear}, such behavior is seen to be possible:

\begin{theorem}\label{theorem - aribtrary growth}
Suppose that $\{c_n\}$ and $\{d_n\}$ are two increasing sequences of positive real numbers, both in $o(n)$, the differences \[\Delta c_n = c_{n+1}-c_n\] are in $O(1)$ (similarly for $\{\Delta d_n\}$), and at least one of $\{c_n\}$, $\{d_n\}$ is divergent.  Then there is a dense set of $\theta$ such that if $\{c_n\}$ is divergent, then
\[\limsup_{n \rightarrow \infty} \frac{M_n(0)}{c_n} = 1,\] while if $\{c_n\}$ is bounded then so is $M_n(0)$.  Similarly, if $\{d_n\}$ is divergent, then
\[\limsup_{n \rightarrow \infty} \frac{|m_n(0)|}{d_n}=1,\] while if $\{d_n\}$ is bounded then so is $m_n(0)$.
\end{theorem}

A closely related result concerns the sequence of values $M_n(x)/|m_n(x)|$:
\begin{theorem} \label{theorem - arbitrary interval for ratio}
Let $0 \leq  r_1 \leq r_2 \leq \infty$.  Then there is a dense set of $\theta$ such that the set of accumulation points of the sequence
\[\left\{ \frac{M_n(0)}{|m_n(0)|}: n=0,1,2,\ldots\right\}\] is the interval $[r_1,r_2]$.
\end{theorem}

We will also include a partial rederivation of \cite[Theorem 1]{MR2515388} in Corollary \ref{corollary - heaviness criterion}: a characterization of those $\theta$ for which $S_n(\theta)\geq 0$ for all $n\geq 0$.

A classical application of the Denjoy-Koksma inequality is that if the $a_i(\theta)$ are drawn from a finite set (such $\theta$ are said to be \textit{badly approximable} or \textit{of finite type}), then $S_n(x) \in O( \log n)$.
\begin{theorem}\label{theorem - bounded type are log}
If $\theta$ is of finite type, then for all $x$ we have $\rho_n(x) \sim \log n$, meaning that the ratio is bounded away from both zero and infinity.
\end{theorem}

\begin{corollary}
If $\theta$ is of finite type, then $|S_n(x)| \notin o(\log n)$ for every $x$, and
\[m_n(x) \in o(\log n) \quad \Longrightarrow \quad M_n(x) \sim \log n,\] and vice-versa.
\end{corollary}

If $A\cup B$ represents a single interval, then as $S^1$ has been partitioned into two intervals of length $\theta$ and $1-\theta$, the analogous problem would be to encode the \textit{Sturmian sequences}, and generating Sturmian sequences using a sequence of substitutions is intimately related to continued fraction expansions for numbers: see for example \cite[Chapter 6]{MR1970385}.  The study of substitutions as they relate to discrepancy sequences of different intervals has been initiated before \cite{MR2040589}, in this paper our approach is different:
\begin{itemize}
\item the interval $[0,1/2]$ is not \textit{dynamically defined}, i.e. not dependent on $\theta$ (although it is fixed),
\item we develop an approach for all $\theta$ (not just quadratic surds, though the process is nicest in this setting),
\item we generate the orbit of any starting point $x$ (though $x=0$ is one particularly nice case that we investigate).
\end{itemize}

\section{Symbol Spaces, Encodings, and Substitutions}\label{section - symbol background}
All background material pertaining to common definitions in symbolic dynamics and substitution systems may be found in \cite[Chapter 1]{MR1970385}; we present here only a short summary of specific notation used herein.  Let $\mathcal{A}=\{A,B,C\}$, and denote by $\mathcal{A}^*$ the \textit{free monoid} on $\mathcal{A}$.  Given $\omega  \in \mathcal{A}^*$, we denote
\[\omega = (\omega)_0 (\omega)_1 \ldots (\omega)_{n-1},\] and say that $\omega$ is a \textit{word of length $n$} with \textit{letters} $(\omega)_i$ drawn from the \textit{alphabet} $\mathcal{A}$.  Note that $\omega_i$ will refer to a sequence of words indexed by $i$, while $(\omega)_i$ will denote the individual letters of a fixed word $\omega$.  This similarity is a potential source of confusion, but the latter notation is much more common in this work: we will rarely refer to specific letters in a given word.

Denote by $|\omega|$ the length of $\omega$.  Elements in $\mathcal{A}^*$ multiply by concatenation, and we adopt power notation for this operation: $(AB)^3=ABABAB$, for example.  The empty word (the identity under concatenation) we denote $\emptyset$.  A \textit{factor} of $\omega$ (of finite or infinite length) is some finite word $\psi$ of length $n$ such that there is some $i$ for which
\[(\psi)_j = (\omega)_{i+j}, \quad j=0,1,\ldots,n-1.\]  If $i=0$ then we say $\psi$ is an \textit{left factor} of $\omega$, and we say $\psi$ is a \textit{right factor} of $\omega$ if $(\psi)_{n-1}=(\omega)_{|\omega|-1}$.  The factor $\psi$ will be called \textit{proper} if $\psi \notin \{\omega, \emptyset\}$.

Any map $\sigma:\mathcal{A} \rightarrow \mathcal{A}^*$ may be extended to a map on $\mathcal{A}^*$ be requiring it to be a homomorphism.  The following is nonstandard but natural.  Endow $\mathcal{A}^{\mathbb{N}}$ with the cylinder topology, and let a finite word $\omega \in \mathcal{A}^*$ represent a clopen set: the set of all elements of $\mathcal{A}^{\mathbb{N}}$ with left factor $\omega$.  We may then further extended $\sigma$ to a map on $\mathcal{A}^{\mathbb{N}}$ by defining
\[\sigma(\omega) = \bigcap_{i=0}^{\infty}\sigma((\omega)_0 (\omega)_1 \ldots (\omega)_{i-1}).\] In all of these situations we refer to $\sigma$ as a \textit{substitution}.

Given a sequence of words $\omega_0,\omega_1,\ldots$ such that $\omega_i$ is a left factor of $\omega_{i+1}$, if
\[ \bigcap_{i=0}^{\infty} \omega_i = \{x\}, \] then we say that $x \in \mathcal{A}^{\mathbb{N}}$ is the \textit{limit} of the words $\omega_i$.

Now consider the space $S^1=[0,1)$ with the map $R_{\theta}(x)=x+\theta \mod 1$ for some irrational $\theta$.  Suppose that $X$ is partitioned into three intervals $A$, $B$, and $C$.  Then given a word $\omega$, we say that $\omega$ \textit{encodes the orbit of $x$} if for all $i \leq |\omega|-1$ we have
\[(\omega)_{i} = A \quad \Longleftrightarrow \quad x+i \theta \in A,\] and similarly for $B$ and $C$.  Given a partition, then, to each $x \in S^1$ we may identify an infinite word $\omega \in \Omega$: the infinite word which encodes the (forward) orbit of $x$.

Let $\mathcal{D}$ be the discontinuities of $(f \circ R_{\theta}^i)(x)$ for $i = 0,1,2,\ldots$:
\[\mathcal{D}=\{-i \theta, -i \theta+ 1/2\}, \quad i=0,1,2,\ldots.\]  For each $x \in \mathcal{D}$, then, we replace $x \in S^1$ with two points, a right and left limit, denoted $x^+$ and $x^-$.  We set
\[R_{\theta}(0^+)=R_{\theta}(1^-)=\theta,\] and similarly for $(1/2)^{\pm}$; while this makes the rotation two-to-one at these points, note that with respect to the alphabet $\mathcal{A}$, the symbolic coding for the forward orbit of $\theta^+$ and $\theta^-$ are identical, so we do not distinguish them.  We still denote our space by $S^1$.  We may now make each of $A$, $B$ and $C$ closed, although we have made $S^1$ totally disconnected.

Given an irrational $\theta$, partition $S^1=[0^+,1^-]$ according to Table \ref{table - parititon} and in a slight abuse of notation let $S^1$ be the set of all words which encode orbits with respect to these conventions.

\begin{table}[bht]
\begin{tabular}{| l | l |}
\hline $\theta<1/2$ & $\theta>1/2$\\
\hline \hline $A=\left[ 0^+, \frac{1}{2}^- \right]$ & $C=\left[ 0^+, (1-\theta)^-\right]$\\
 $B=\left[\frac{1}{2}^+,(1-\theta)^- \right]$ & $B=\left[(1-\theta)^+,\frac{1}{2}^- \right]$\\
$C = \left[(1-\theta)^+,1 \right]$ & $A=\left[\frac{1}{2}^+,1^- \right]$\\
\hline
\end{tabular}
\caption{The partition $S^1=A \cup B \cup C$ depending on $\theta$.}
\label{table - parititon}
\end{table}

The following lemma is immediate, and immediately explains the apparent ambiguity in the statement of Theorem \ref{theorem - substitution sequence}:
\begin{lemma}\label{lemma - at most two errors}
If $\omega$ is an infinite word encoding the orbit of a point $x \in S^1$ under rotation by $\theta$, then $\omega$ encodes the orbit of some $x \in S^1$ \textit{without} the introduction of $\mathcal{D}$ with at most two errors.  Alternately the coding is exact up a change of endpoints of the intervals $A$, $B$ and $C$.
\begin{proof}
The orbit of any point can hit the endpoints of $A$, $B$ and $C$ at most twice.
\end{proof}
\end{lemma}

\section{The Renormalization Procedure}
Recall $\gamma$, the Gauss map \eqref{eqn - gauss map}; we define a similar map.
\begin{equation}\label{eqn - g}
g([a_1,a_2,a_3,\ldots]) = \begin{cases}[a_3,a_4,\ldots]=\gamma^2(\theta) & (a_1=0 \bmod 2)\\ [1,a_2,a_3,\ldots]=\frac{1}{1+\gamma(\theta)} & (a_1 = 1 \bmod 2, \, a_1 \neq 1)\\ [a_2+1,a_3,\ldots]=1-\theta & (a_1=1).\end{cases}
\end{equation}
Note that if $\theta>1/2$, then necessarily $g(\theta)<1/2$.  It will be convenient to define
\begin{equation}\label{eqn - E}
E(x) = \max\{n \leq x: n \in \Z, \, n=0 \bmod 2\}.\end{equation}

The triplet $\{X, \mu, T\}$ refers to a compact probability space $\{X, \mu\}$ and a continuous transformation $T$ on $X$ which preserves $\mu$.  Given irrational $\theta$, we denote
\begin{equation}\label{eqn - stuff related to theta_n}
\theta_n  = g^n(\theta),\quad \delta_n  = 1-E(a_1(\theta_n))\theta_n,\quad I_n = \{S^1, \mu, R_{\theta_n}\}.
\end{equation}
Note that $\delta_n=1$ if and only if $\theta>1/2$; otherwise $\delta_n<1/2$.

Partition each $I_n$ into intervals $A$, $B$ and $C$ according to Table \ref{table - parititon}, and recall that by convention we have disconnected each $I_n$ such that all iterates of the characteristic functions of $A$, $B$ and $C$ under $R^i_{\theta_n}$ are continuous.  Given $\{X, \mu, T\}$ and a set $S \subset X$, the \textit{return time to $S$} is given by
\[n(x) = \min\{n>0: T^n(x) \in S\}.\]
As irrational rotations are minimal, $n(x)$ will be defined for all $x \in S^1$ if $S$ is an interval of positive length.  The \textit{induced system} on $S$ is defined by \[\{S, \mu|_S, T|_s\},\] where $T|_S (x) = T^{n(x)}(x)$ for all $x \in S$.  Define $I_{n+1}' \subset I_n$ by
\[I'_{n+1}=[0^+, \delta_n^-].\]

Finally, define the substitutions $\sigma_n=\sigma(\theta_n)$ according to Table \ref{table - substitutions}, and define the functions $\varphi_n=\varphi(\theta_n)$ according to:
\begin{equation}\label{eqn - varphi}
\varphi(x) = \begin{cases} 1- x & (a_1(\theta)=1)\\ \delta_n^{-1}x & (a_1(\theta) \neq 1)\end{cases}
\end{equation}

\begin{table}[hbt]
\begin{tabular}{| l | l |}
\hline Case & Substitution\\
\hline \hline \multirow{3}{*}{$a_1=2k, \, a_3 \neq 1$} &$A \rightarrow (A^{k+1}B^{k-1}C) (A^{k}B^{k-1}C)^{a_2-1}$\\
& $ B \rightarrow (A^{k}B^{k}C) (A^{k}B^{k-1}C)^{a_2-1}$\\
& $C \rightarrow (A^{k}B^{k}C) (A^{k}B^{k-1}C)^{a_2}$\\
\hline \multirow{3}{*}{$a_1=2k, \, a_3=1$} & $A \rightarrow (A^{k}B^{k}C) (A^{k}B^{k-1}C)^{a_2}$\\
& $B \rightarrow (A^{k+1}B^{k-1}C)(A^{k}B^{k-1}C)^{a_2}$\\
& $C \rightarrow (A^{k+1}B^{k-1}C)(A^{k}B^{k-1}C)^{a_2-1}$\\
\hline \multirow{3}{*}{$a_1=2k+1$}& $A \rightarrow A^{k}B^{k}C$\\
&$B \rightarrow A^{k+1}B^{k-1}C$\\
&$C \rightarrow A$\\
\hline \multirow{3}{*}{$a_1=1$} & $A \rightarrow A$\\
&$ B \rightarrow B$\\
&$C \rightarrow C$\\
\hline
\end{tabular}
\caption{The substitution $\sigma$ as a function of $\theta$.}
\label{table - substitutions}
\end{table}

\begin{lemma}\label{lemma - balanced word}
Suppose that $\theta < 1/2$, $E(a_1(\theta))=2k$, and
\[(1-2k\theta)^+ \leq x \leq \left(\frac{1}{2}-(k-1)\theta\right)^-.\]
Then the orbit of $x$ begins $A^k B^{k-1} C$.
\begin{proof}
The assumption $\theta<1/2$ tells us how to partition $S^1$ according to Table \ref{table - parititon} as well as guaranteeing that $k \geq 1$. Note that the lower inequality certainly guarantees that
\[\frac{1}{2}-k \theta < x \leq \left(\frac{1}{2}-(k-1)\theta\right)^-,\] which tells us that $x+i \theta \leq (1/2)^-$ for $i=0,1,\ldots (k-1)$, while $x+k\theta>1/2$.  So the coding of the orbit of $x$ begins with exactly $A^k$ before seeing either $B$ or $C$.  As we know
\[\left(1-2k\theta\right)^+\leq x < 1- (2k-1)\theta,\] we know that we have $x+(2k-1)\theta<1$, while $x+2k\theta\geq 1^+$.  Therefore, once we have accounted for the points
$x+i\theta$ for $i=0,1,\ldots,k-1$, the terms $i=k,k+1,\ldots,(2k-1)$ must all belong to either $B$ or $C$.  That $C$ is an interval of length exactly $\theta$ guarantees that exactly the final term is $C$.  The rest of the terms (if there are any) are therefore $B$.
\end{proof}
\end{lemma}

\begin{proposition}\label{proposition - setup for induction}
We have the measurable and continuous isomorphism
\[\left\{ I'_{n+1}, \mu|_{I'_{n+1}}, \left(R_{\theta_n}\right)|_{I_{n+1}}\right\} \xrightarrow{\varphi_n} \left\{I_{n+1}, \mu, R_{\theta_{n+1}}\right\}.\]  Furthermore, for all $x \in A \subset I_{n+1}$, the word $\sigma_n(A)$ encodes the orbit of $\varphi^{-1}(x)$ through its return to $I_{n+1}'$ (the encoding is with respect to the partition $A$, $B$, $C$ in $I_n$), and similarly for $B$ and $C$.
\begin{proof}
In the case that $\theta_n>1/2$, then $\theta_{n+1}=1-\theta_n$ and $I_{n+1}'=[0^+,1^-]$.  However, by referring to Table \ref{table - parititon}, we see that the intervals $A$, $B$ and $C$ exactly reflect the reversal of orientation given by $\varphi_n(x)=1-x$, and the substitution $\sigma_n$ is identity.  So we proceed on the assumption that $\theta_n<1/2$: in $I_n$ we have
\[A=[0^+,1/2^-], \quad B=[1/2^+,(1-\theta)^-], \quad C=[(1-\theta)^+,1^-].\]
Then $\varphi_n$ is scalar multiplication by $\delta_n^{-1}$, so there are only two things to show:
\begin{itemize}
\item The first-return map $(R_{\theta_n})|_{I'_{n+1}}$ is rotation by $\theta_{n+1}$, after rescaling by $\varphi_n$, and
\item the substitution $\sigma_n$ encodes the correct information.
\end{itemize}
There are three cases to consider: $a_1(\theta_n)=1 \bmod 2$, or $a_1(\theta_n)=0 \bmod 2$ with the sub-cases $a_3(\theta_n)=1$ or $\neq 1$.  Assume for now that $a_1(\theta_n)=0\bmod 2 $ and $a_3(\theta_n)=1$.

As $a_1(\theta_n)=0\bmod2$ and $a_3(\theta_n)=1$, we have $g(\theta_n)=\gamma^2(\theta_n)>1/2$, so in $I_{n+1}$ we have
\[C=[0^+,(1-\theta_{n+1})^-], \quad B=[(1-\theta_{n+1})^+,1/2^-], \quad A=[1/2^+,1^-],\] with corresponding preimages in $I'_{n+1}$ scaled by $\delta_n$.  We will first verify that the intervals have the desired return times (which may be read from the length of the words $\sigma_n(A)$, $\sigma_n(B)$ and $\sigma_n(C)$) and that the induced map is indeed rotation by $\theta_{n+1}$ (up to scale $\delta_n$).  As $E(a_1(\theta_n))=a_1(\theta_n)$ we have
\[\delta_n=\|q_1(\theta_n) \cdot \theta_n\|,\] from which it follows that the return time of $0$ is
\[n(0)=q_2=a_1 a_2 +1,\] and one may now verify that the entire interval $\varphi_n^{-1}(C)$ has this return time; the preimage of the right endpoint of $C$ under $\varphi_n$ is exactly $1-(q_1+q_2)\theta_n$.  The remaining points in $I'_{n+1}$ have return time $q_2+q_1$ and the induced map is a rotation by $q_2 \theta_n$ on $[0^+,\delta_n^-]$; see Figure \ref{figure-induced map}.

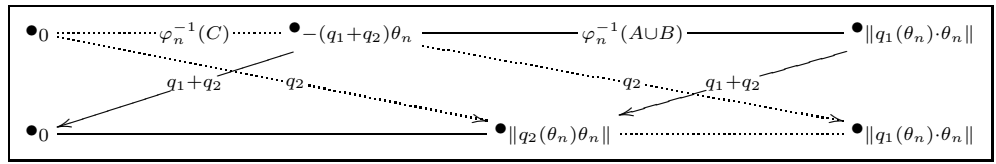
\begin{figure}[h]
\center{\framebox{\xymatrix{
\bullet_0 \ar@{.>}[drrrr]|{\, q_2 \,} \ar@{.}[rrr]|{\, \varphi_n^{-1}(C) \,}& & & \bullet_{-(q_1+q_2)\theta_n}\ar@{-}[rrrr]|{\, \varphi_n^{-1}(A \cup B) \,}\ar[dlll]|{\, q_1+q_2\, } \ar@{.>}[drrrr]|{\, q_2\, } & & & & \bullet_{\|q_1(\theta_n)\cdot\theta_n\|}\ar[dlll]|{\, q_1+q_2 \,} \\
\bullet_0 \ar@{-}[rrrr] & & & &\bullet_{\|q_2(\theta_n) \theta_n\| }\ar@{.}[rrr] & & & \bullet_{\|q_1(\theta_n)\cdot\theta_n\| }
}}}\caption{\label{figure-induced map} \small{\textit{Return times for the case $a_1(\theta_n)=0\bmod 2$, $a_3(\theta_n)=1$.}}}
\end{figure}

At this point we may verify that the rotation is by $g(\theta_n)$, up to scale:
\begin{align*}
\frac{\|q_2(\theta_n)\cdot \theta_n\|}{\delta_n} & = \frac{q_2(\theta_n)\cdot \theta_n -p_2(\theta_n)}{1-q_1(\theta_n)\cdot \theta_n}\\
&= \frac{(a_1 a_2+1)\theta_n-a_2}{1-a_1 \theta_n}\\
&= \frac{a_2\left(a_1-\frac{1}{\theta_n}\right) +1}{\frac{1}{\theta_n}-a_1}\\
&= \frac{1-a_2 \gamma(\theta_n)}{\gamma(\theta_n)}\\
&=\gamma^2(\theta_n).
\end{align*}

Now suppose that $x \in \varphi^{-1}(B)$, and for convenience denote $E(a_1)=a_1=2k$.  Clearly, the orbit of $x$ begins with a point in $A$ (in $I_n$, as $A=[0^+,1/2^-]$ contains $[0^+,\delta_n^-]$).  As $x < 1/2-k \theta_n$, however, we have
\[(1-2k\theta_n)^+ \leq x+\theta_n \leq \left((1/2)-(k-1)\theta_n\right)^-,\] so by Lemma \ref{lemma - balanced word}, we may concatenate the word $A^k B^{k-1} C$ to this initial $A$.  Since $2k=a_1$, we now have
\[x+\theta_n+(2k \theta_n) < x+\theta_n \leq \left( (1/2) - (k-1)\theta_n\right)^-.\]  Either we have returned to $I'_{n+1}$, in which case we are done, or we have not, in which case we apply Lemma \ref{lemma - balanced word} again, repeating until we return to $I'_{n+1}$, which must take a total of $q_2+q_1=a_1 (a_2+1) + 1 $ steps.

For those points in the interval $\varphi_n^{-1}(a)$, note that the only discontinuity of $R_{\theta_n}^i$ for $i=0,1,\ldots,q_2$ to distinguish the orbits compared to points in $\varphi_n^{-1}(A)$ is the point $1/2-k \theta$, which will change the single term $x+k \theta$ from an `$A$' to a `$B$'. Points in $\varphi_n^{-1}(C)$ are considered identically to those in $\varphi_n^{-1}(B)$, noting that the shorter return time requires one fewer concatenation of $A^{a_1}B^{a_1-1}C$.

The other cases are similarly considered; the case $a_1(\theta_n)=0\bmod 2$, $a_3(\theta_n) \neq 1$ is nearly identical, while for the case $a_1(\theta_n)=1 \bmod 1, \, \neq 1$ we have $\delta_n > \theta_n$, so the return time of $0^+$ is one, explaining the much shorter substitution $\sigma_n(C)=A$ in this case.
\end{proof}
\end{proposition}
Denote the iterated pull-back of $I_n$ into $I_0$ by
\begin{equation}\label{eqn - tilde I}\tilde{I}_n = \left(\varphi_0^{-1} \circ \cdots \circ \varphi_{n-1}^{-1}\right)(I_n).\end{equation}
\begin{corollary}\label{corollary - iterated substitutions}
We have the measurable and continuous isomorphism
\[\left\{ \tilde{I}_n, \mu|_{\tilde{I}_n}, \left(R_{\theta}\right)|_{\tilde{I}_n}\right\} \xrightarrow{(\varphi_{n-1} \circ \cdots \circ \varphi_0)} \left\{I_n, \mu, R_{\theta_n}\right\}.\]
Furthermore, for any $x \in A \subset I_n$, the word $\left(\sigma_0 \circ \cdots \circ \sigma_{n-1}\right)(A)$ encodes the orbit of $\left(\varphi_{0}^{-1}\circ \cdots \circ \varphi_{n-1}^{-1}\right)(x)$ in $I_0$ through its return to $\tilde{I}_n$, and similarly for $B$, $C$.
\end{corollary}
\section{Proof of Theorem \ref{theorem - substitution sequence}}

The proof of \eqref{eqn - special substitution point} is immediate in light of Corollary \ref{corollary - iterated substitutions}; the point $x(\theta)$ is given by
\[x(\theta) = \bigcap_{i=0}^{\infty} \tilde{I}_i,\] where the $\tilde{I}_i$ were defined in \eqref{eqn - tilde I}.  This intersection is nonempty as the sets are nested closed intervals in the compact space $S^1$.  The length of $\tilde{I}_n$ is given by
\[\delta_0 \cdot \delta_1 \cdots \delta_{n-1},\] and we have already remarked that for $\theta_n<1/2$, we have $\delta_n < 1/2$.  As no two successive terms in the sequence $\theta_0, \theta_1, \ldots$ may be larger than one half, the length tends to zero, and the intersection is either a singleton or a pair $\{x^-,x^+\}$.  In the latter scenario, however, both $x^-$ and $x^+$ would have identical coding of their forward orbits.  As we did not `split' the points $i \theta$ or $i\theta+1/2$ for $i>0$ when disconnecting $S^1$, this is not possible.

As all non-identity substitutions map each letter to a word beginning in $A$, and all non-identity substitutions map $A$ to a word of length at least three, and no two consecutive substitutions may be identity, it follows that the sequence of words
\[\left( \sigma_0 \circ \sigma_1 \circ \cdots \circ \sigma_{n-1}\right)(\omega)\] has a limit regardless of the choice of nonempty $\omega$, and Corollary \ref{corollary - iterated substitutions} shows that this word must encode the orbit of $x(\theta)$ in the disconnected version of $S^1$.  Lemma \ref{lemma - at most two errors} finishes the proof of this portion of Theorem \ref{theorem - substitution sequence}.

Let us now turn our attention to constructing the orbit of an arbitrary $x_0 \in S^1$.  Define
\[x_1 = x_0+i \theta, \quad i \in \{j \geq 0 : x+j \theta \in I'_1\},\] and let $\omega_0$ be the word which encodes the orbit of $x_0$ through its arrival to $x_1$; if $x_0 \in I'_1$, we may set $\omega_0$ to be the empty word (though we are not required to do so).  We now pass to the system $I_1$, letting $(x_1 \in I_1) = \varphi_0(x_1 \in I'_1)$.  We set $x_2$ to be a point in $I'_2$ which is in the orbit of $x_1$, and let $\omega_1$ be the word encoding this finite portion of the orbit, then pass to $I_2$, etc.  Equation \eqref{eqn - arbitrary substitution point} now follows from Proposition \ref{proposition - setup for induction} so long as infinitely many $\omega_n \neq \emptyset$.  We only have the option of letting all but finitely many $\omega_n$ be empty if $x$ is a preimage of $x(\theta)$; we have already remarked in this case that the limiting word may be found handily.

A potential source of confusion at this point is the desire to claim that $x(\theta)=0$, as we always construct $I'_{n+1}=[0^+,\delta_n^-]$.  However, $\varphi_n(x)=1-x$ for those $n$ such that $\theta_n>1/2$.  So $\varphi_n^{-1} \circ \varphi_{n+1}^{-1}$ pulls back $I_{n+2}$ to the interval $[(1-\delta_{n+1})^+,1^-] \subset I_n$.  Those $\theta$ for which $x(\theta)=0$ will be addressed in Proposition \ref{proposition - heavy theta have zero}.

\begin{proposition}\label{proposition - canonical choice for omega_n}
Without loss of generality, $\omega_n$ may be required to either be empty, or a proper right factor of either $\sigma_n(A)$, $\sigma_n(B)$, or $\sigma_n(C)$.
\begin{proof}
The images of $R_{\theta_n}^i\left(I'_{n+1}\right)$ cover all of $I_n$ through the return times, so any $x$ may be viewed as returning to $I'_{n+1}$ via a right factor of one of these words.  If the return is through the entire word $\sigma_n(A)$, we would have begun with $x_n \in I'_{n+1}$ and could have set $\omega_n = \emptyset$.
\end{proof}
\end{proposition}
\begin{remark}
One could alternately require that $\omega_n$ be nonempty by allowing all nonempty right factors of $\sigma_n(A)$, $\sigma_n(B)$, and $\sigma_n(C)$; instead of $\omega_n=\emptyset$ for $x \in I'_{n+1}$, let $\omega_n$ be $\sigma$ applied to the letter encoding whichever interval in $I_{n+1}$ contains $\varphi_n(x)$.
\end{remark}

In order to construct the orbit of zero we will side-step this computation altogether:

\begin{lemma}\label{lemma - trick for zero}
Suppose that $\theta_n>1/2$.  Let $\Omega$ encode the orbit of $0^+$ in the system $I_n$, and $\Upsilon$ encode the orbit of $0^+$ in the system $I_{n+1}$.  Then for all $i \geq 1$, $(\Omega)_i=(\Upsilon)_i$.  For $i=0$, $(\Omega)_0=C$ while $(\Upsilon)_0=A$.
\begin{proof}
The isomorphism $\varphi_n(x)=1-x$ and the identity substitution $\sigma_n$ ensures that $\Omega$ is identical to the coding of the orbit of $1^-$ in $I_{n+1}$.  As the forward orbit of $0$ under rotation by the irrational $\theta_n$ does not hit any other endpoints of the intervals $A$, $B$, and $C$, we have that the orbit of $1^-$ and $0^+$ in the system $I_{n+1}$ are identical after this initial term.
\end{proof}
\end{lemma}

With this lemma in mind, then, define the map $\Psi(\omega)$ on both $\mathcal{A}^*$ and $\mathcal{A}^{\mathbb{N}}$:
\begin{equation}\label{eqn - psi}
\left(\Psi \omega\right)_i = \begin{cases} C & (i=0)\\ \omega_i & (i \neq 0).\end{cases}\end{equation}

Define the maps $\sigma'_n=\sigma'(\theta_n)$:
\begin{equation}\label{eqn - sigma prime}
\sigma'(\theta) = \begin{cases} \sigma(\theta) & (\theta<1/2)\\ \Psi & (\theta>1/2).\end{cases} \end{equation}
Then \eqref{eqn - zero substitutions} follows if we appropriately choose the words $\omega'_n$ to accurately encode some string of the initial orbit of $0^+$ in $I_n$.  Then the resulting word
\[\left(\sigma_0' \circ \sigma_1' \circ \cdots \circ \sigma_{n-1}'\right)(\omega'_n)\] will accurately represent the initial orbit of $0^+$, but it is no longer guaranteed that the length of this word increases!  For example, if $\theta=[3,2,2,2,2,\ldots]$, then we will alternate between $\sigma_n'$ being $\Psi$ and a substitution which maps $C \rightarrow A$.  Setting $\omega'_n=A$ for all those $n$ for which $\theta_n<1/2$ would therefore always map via this long string of compositions to
\[A \xrightarrow{\Psi} C \xrightarrow{\sigma} A \xrightarrow{\Psi} C \xrightarrow{\sigma} \cdots\]

Define
\begin{equation}\label{eqn - omega' words for zero orbit}
\omega'_n = \begin{cases} A^{k+1} B^{k-1} C & (a_1(\theta_n)=2k)\\ A^{k+1}B^k & (a_1(\theta)=2k+1) \\ \Psi(\omega'_{n+1}) & (a_1(\theta)=1). \end{cases}
\end{equation}

The reader may verify that the word $\omega'_n$ does accurately encode some initial portion of the orbit of $0^+$ depending on the parity of $a_1(\theta_n)$.  Note that whenever $\Psi$ is applied, it affects only the first letter of its input.  From this it follows that if $\omega = (\omega)_0 \nu$, then
\begin{equation}\label{eqn - sigma' pull off one letter}
\left( \sigma'_0 \circ \cdots \circ \sigma'_{n-1}\right)(\omega) = \left( \sigma'_0 \circ \cdots \circ \sigma'_{n-1}\right)((\omega)_0) \left( \sigma_0 \circ \cdots \circ \sigma_{n-1}\right)(\nu).
\end{equation}

As $\omega'_n$ always has length larger than one, our previous reasoning now guarantees that the length of $\Omega'_n$ diverges, establishing \eqref{eqn - zero substitutions} and completing the proof. \qedhere

Before moving on to the study of the growth rates of discrepancy sums, we present a few observations about this process.

\begin{proposition}\label{proposition - heavy theta have zero}
Those $\theta$ for which $x(\theta)=0 (=0^+)$ are exactly the set
\begin{equation}\label{eqn - heavy theta}H = \left\{ \theta: a_{2i-1}(\theta)=0 \bmod 2, \, i=1,2,\ldots\right\}.\end{equation}
\begin{proof}
We leave the reader to verify that $H$ is exactly the set of $\theta$ for which $g^n(\theta)<1/2$ for every $n$.  For those $\theta \in H$, then, we always have $I'_{n+1}=[0^+,\delta_n^-]$, where $\delta_n < 1$, and we never need apply the isomorphism $\varphi_n(x)=1-x$.  That is,
\[0 \in \left( \varphi_0^{-1} \circ \cdots \circ \varphi_{n-1}^{-1}\right)(I_n)\] for all $n$: $0 = x(\theta)$.

On the other hand, if $n$ is the first index such that $\theta_n > 1/2$, we must have $\varphi_n(x)=1-x$.  As $\theta_{n+1}<1/2$, however, it follows that within $I_n$, we have
\[\varphi_{n}^{-1} \circ \varphi_{n+1}^{-1} (I_{n+2}) = [(1-\delta_{n+1})^+,1^-],\] from which it follows that
\[0 \notin \left( \varphi_0^{-1} \circ \cdots \circ \varphi_{n+1}^{-1}\right)(I_{n+2}). \qedhere \]
\end{proof}
\end{proposition}

\begin{proposition}\label{proposition - quadratics are periodic}
The sequence of substitutions $\sigma_n$ is eventually periodic if and only if $\theta$ is a quadratic surd.
\begin{proof}
Clearly the sequence $\sigma_n$ is eventually periodic if and only if the orbit of $\theta$ under $g$ is eventually periodic.  From the definition \eqref{eqn - g} of $g$ we have for all $i \geq 2$
\begin{equation}\label{eqn - g shifts terms}a_i(\theta_{n+1}) = a_{i+k}(\theta_n): \quad k = \begin{cases} 0 & (a_1(\theta_n)=1\bmod 2, \,\neq 1)\\ 1 & (a_1(\theta_n)=1)\\ 2 & \left(a_1(\theta_n)=0\bmod 2\right)\end{cases}\end{equation}
So, if $a_i(\theta)$ are eventually periodic (Gauss' criteria for quadratic surds), we must have infinitely many $n$ such that for all $i \geq 2$ we have for any $j,k$ \[a_i(\theta_{n_k})=a_i(\theta_{n_j}).\]
Suppose that a period of $a_i(\theta)$ is given by the terms $\alpha_1, \ldots, \alpha_N$, and assume without loss of generality that for $i \geq 2$
\[a_i(\theta_{n_k})=\alpha_{i \bmod N}.\]  Then $a_1(\theta_{n_k})$ is either $1$, $\alpha_1$, or $\alpha_1+1$.  Since the collection $n_k$ was infinite, one value must be taken twice, giving a period in the orbit $g(\theta)$.

On the other hand, assume that $\theta_j=\theta_{j+nk}$ for $n=0,1,\ldots$ and $k \neq 0$.  From \eqref{eqn - g shifts terms} it follows that $a_i(\theta)$ is eventually periodic.
\end{proof}
\end{proposition}

\begin{remark}
The periods under $g$ and $\gamma$ need not be the same, nor is one necessarily longer than the other.  For example, the golden mean has period one under $\gamma$ but period two under $g$, while $\theta=[2,1,2,1,\ldots]$ has period two under $\gamma$ and period one under $g$.  Furthermore, the sequence $\sigma_n$ is purely periodic if and only if $\theta_n=\theta_0$ for some $n \neq 0$, which is not the same as the partial quotients of $\theta$ being purely periodic.  Consider for example $\theta=[3,2,2,2,\ldots]$, whose partial quotients are clearly not purely periodic, but satisfies $\theta_2=\theta_0$.
\end{remark}

\section{The Arithmetic of Our Substitutions}

Let $\theta_0<1/2$, so that
\[f(x) = \begin{cases} +1 & (x \in A)\\ -1 & (x \in B \cup C).\end{cases}\]  For $\theta_0>1/2$ we could repeat all future arguments with a sign change.  Given $\omega \in \mathcal{A}^n$, define (consistent with existing notation)
\begin{align*}
S(\omega) &= \sum_{i=0}^{n-1} \left( \chi_A - \chi_{B \cup C}\right)\omega_i,\\
M(\omega) &= \max \left\{S(\omega_0 \ldots \omega_{j-1}): j=1,2,\ldots,n\right\},\\
m(\omega) &= \min \left\{S(\omega_0 \ldots \omega_{j-1}): j=1,2,\ldots,n\right\}.
\end{align*}
Note that we \textit{do not} include the empty word in determining $M(\omega)$, $m(\omega)$.
\begin{proposition}\label{proposition - arithmetic one step}
Suppose $|\omega|=n \neq 0$, $\omega \neq C$, $M(\omega) \geq 0$, $\omega$ does not have $CC$, $CB$ or $BA$ as factors, and $\sigma$ is a substitution given by Table \ref{table - substitutions}, depending on $\theta$.  If $a_1(\theta)=0 \bmod 2$ and $a_3(\theta) \neq 1$, or if $a_1(\theta)=1$, then:
\[S(\sigma(\omega))=S(\omega), \quad M(\sigma(\omega))=M(\omega)+E(a_1), \quad m(\sigma(\omega))=m(\omega).\]

On the other hand, if $a_1(\theta)=0 \bmod 2 $ and $a_3(\theta)=1$, then
\[S(\sigma(\omega))=-S(\omega), \quad M(\sigma(\omega))=-m(\omega)+E(a_1), \quad m(\sigma(\omega))=-M(\omega).\]
Finally, if $a_1(\theta)=1 \bmod 2$, $\neq 1$, and either
\begin{itemize}
\item $(\omega)_{n-1}\neq C$, or
\item $(\omega)_{n-1}=C$, but there is some $j \neq n$ such that $S((\omega)_0 (\omega)_1 \ldots (\omega)_{j-1})=m(\omega)$,
\end{itemize}
then also
\[S(\sigma(\omega))=-S(\omega), \quad M(\sigma(\omega))=-m(\omega)+E(a_1), \quad m(\sigma(\omega))=-M(\omega).\]
If $a_1(\theta)=1 \bmod 2$, $(\omega)_{n-1}=C$ and $S((\omega)_0 \ldots (\omega)_{j-1})>m(\omega)$ for all $j \neq n$, then
\[S(\sigma(\omega))=-S(\omega), \quad M(\sigma(\omega))=-m(\omega)-1+E(a_1), \quad m(\sigma(\omega))=-M(\omega).\]
\begin{proof}
The prohibition on $CB$, $CC$ and $BA$ being factors of $\omega$ are necessary for $\omega$ to encode the orbit of any point under rotation by any $\theta$, so this condition is not prohibitive in our setting.

In all cases, the statements regarding the value $S(\sigma(\omega))$ follow from examining $S(\sigma(x))$ for each $x \in \mathcal{A}$; the reader may consult Table \ref{table - substitutions} to verify that $S(\sigma(x))=\pm S(x)$ as described, and the statement then follows from the fact that $\sigma$ is a homomorphism.  We will turn our attention, then, to the statements regarding $m(\sigma(\omega))$ and $M(\sigma(\omega))$.  All cases but the last are considered similarly with the possible sign-change outlined above in mind.

For example, suppose that $a_1=0 \bmod 2$ and $a_3 \neq 1$.  Let $\omega=\upsilon \psi$, where $\upsilon$ is the largest left factor of $\omega$ such that $S(\upsilon)=M(\omega)-1$: note that as $M(\omega) \geq 0$ and the empty word was not considered in computation of $M(\omega)$, we have $(\psi)_0=A$.  As $S(\sigma(\upsilon))=S(\upsilon)=M(\omega)-1$ and $M(\sigma(A))=E(a_1)+1$, we know that
\[M(\sigma(\omega)) \geq M(\sigma(\upsilon)\psi) = M(\omega)+E(a_1).\]

Assume on the other hand that
\[\sigma(\omega)=\sigma(\upsilon)\nu \psi, \quad S(\sigma(\upsilon)\nu)>M(\omega)+E(a_1),\] and $\upsilon$ is of maximal length to allow such a decomposition.  Note that $\nu \neq \emptyset$ as $S(\sigma(\upsilon))=S(\upsilon)\leq M(\omega)$.  As $\upsilon$ is a proper factor, it is followed by a letter, and by maximality on the length of $\upsilon$, $\nu$ is a proper left factor of either $\sigma(A)$, $\sigma(B)$, or $\sigma(C)$, and $E(a_1)\neq 0$.  If $\upsilon$ is followed by $A$ in $\omega$,
\[S(\sigma(\upsilon)) =S(\upsilon) \leq M-1.\]  On the other hand, $S(\nu)\leq E(a_1)+1 = M(\sigma(A))$, contradicting the value $S(\sigma(\upsilon)\nu)$.  The possibility of $\upsilon$ followed by $B$ or $C$ are similarly considered; the larger possible $S(\sigma(\upsilon))=M(\omega)$ is countered by $S(\nu)\leq E(a_1)$ in these cases.

The ambiguity in the situation when $a_1(\theta)=1 \bmod 2$, $\neq 1$ is due to the substitution $\sigma(A)=C$, which does not achieve an intermediate sum of $E(a_1)$ (as does $\sigma(B)$).  On the assumption that there is some proper left factor $\psi$ of $\omega$ such that $S(\psi)=m(\omega)$, however, we know that the letter which follows $\psi$ must be $A$; similar computations to the above then apply.  If the only left factor of $\omega$ which achieves a sum of $m(\omega)$ is in fact $\omega$ itself, then if the final letter of $\omega$ is $B$ we again have no problem.

Assume, then, that $S(\omega)=m(\omega)$, there is no proper left factor with this sum, and $\omega$ ends with the letter $C$.  As $M(\omega) \geq 0$ by assumption, there is a letter preceding this terminal $C$ (that is, $\omega \neq C$).  If this letter is $A$, then the left factor $\psi$ such that $\omega=\psi A C$ has the minimal sum as its sum (even if it is empty), and the preceding reasoning applies.  Therefore $\omega$ must be of the form $\psi BC$ (recall that $CC$ is not a factor): considering $\sigma(B)$ following $S(\sigma(\psi))=-m(\omega)-2$ completes the proposition.
\end{proof}
\end{proposition}
For convenience, denote
\begin{eqnarray}\label{eqn - iterated substitution notation}\sigma^{(n)} &= \sigma_0 \circ \sigma_1 \circ \cdots \circ \sigma_{n-1},\\ \sigma'^{(n)} &= \sigma'_0 \circ \sigma'_1 \circ \cdots \circ \sigma'_{n-1}.\end{eqnarray}

Recall \eqref{eqn - omega' words for zero orbit} and define for $n \geq 1$
\begin{equation}\label{eqn - big Omega words}
\Omega_n = \sigma^{(n)}(A), \quad \Omega'_n = \sigma'^{(n)}(\omega'(n)).
\end{equation}

Define $p_n$ to track the parity of how many $\theta_i>1/2$:
\begin{equation}\label{eqn - p_n}
p_n = \left(\sum_{i=1}^{n-1} \chi_{(1/2,1)}(\theta_i)\right) \bmod 2.
\end{equation}
We now have all the tools necessary to precisely study the sequences $M_n(y)$ and $m_n(y)$ for $y \in \{x(\theta),0\}$:
\begin{proposition}\label{proposition - arithmetic on substitutions}
Assume that $\theta_0<1/2$.  Then
\[S(\Omega_n)  = (-1)^{p_n}, \quad  S(\Omega'_n) = 1\]
\begin{gather*}
\left| M(\Omega_n) - \left( 1 + \sum_{\substack{i\leq n-1 \\ p_i=0}} E(a_1(\theta_i))\right)\right| \leq 1,\quad  M(\Omega'_n) = 1 + \sum_{\substack{i \leq n \\ p_i=0}}E(a_1(\theta_i)),\\
\left| m(\Omega_n) - \left( 1 - \sum_{\substack{i\leq n-1 \\ p_i=1}} E(a_1(\theta_i))\right)\right| \leq 1, \quad m(\Omega'_n) = 1 - \sum_{\substack{i \leq n \\ p_i=1}}E(a_1(\theta_i)).
\end{gather*}
\begin{proof}
The word $\Omega_n$ in \eqref{eqn - big Omega words} is formed by successive substitutions acting on the word $A$; as such, it will always begin with $A$, so $M(\Omega_n)\geq 1$.  We immediately see that all $S(\Omega_n)=\pm 1$ according to the parity of $p_n$ by applying Proposition \ref{proposition - arithmetic one step} in succession.  The ambiguous case in Proposition \ref{proposition - arithmetic one step} arose when $\omega$ was a word which had a nonnegative maximal sum (as do all $\Omega_n$) and whose minimum sum is only achieved as its total sum, with $C$ as a terminal factor.  Furthermore, we would need $\theta_n$ to have first partial quotient odd and larger than one.  For this to happen with the restriction that all $S(\Omega_n)= \pm 1$ requires that $S(\Omega_n)=-1$ (otherwise the minimal sum is achieved by the proper left factor $A$), and therefore $S(\Omega_{n-1})=1$.  This scenario \textit{also} require that $M(\Omega_{n-1})=1$ (otherwise $m(\Omega_n)<-1 \leq S(\Omega_n)$); so this situation can only occur in our scenario when $\Omega_{n-1}=A$: \textit{this possible error of one may only appear once in the sequence of arithmetic computations from repeated application of Proposition \ref{proposition - arithmetic one step}}.

We leave to the reader the verification that the parity of $p_n$ exactly dictates whether substitutions will add to the maximal values or subtract from the minimal values; refer to Proposition \ref{proposition - arithmetic one step} again.

Let us now consider $\Omega'_n$.  Note that $\sigma_j'=\Psi$ exactly when $\theta_j>1/2$, exactly when $\sigma_{j-1}$ has the property that $S(\sigma_{j-1}(\omega))=-S(\omega)$.  Clearly we have $S(\Psi(\omega))=S(\omega)-2$ provided $\omega$ begins with $A$.  Also note that if $S(\omega)=1$, then if $m(\omega)=1$ we must have $\omega_0=A$: \textit{it is never possible in our construction for $\omega$ to terminate with $C$, $S(\omega)=1$, and $m(\Psi(\omega))=S(\Psi(\omega))$ is the only time this value is reached.}

Our choice of $\omega'(n)$ always begins with $A$ and has $S(\omega'(n))=1$, and for those $\sigma_n$ such that $S(\sigma_n(A))=-1$, the reader may verify that
\[S\left(\sigma_n(\Psi(\omega))\right)=2-S(\omega)\] by applying Proposition \ref{proposition - arithmetic one step}.  While this change will change the sum of $+1$ to $-1$, it is immediately followed by a substitution which reverses the sign of the sum: we maintain
\[S(\Omega'_n)=1.\]

Furthermore, as $m(\omega'_n)=1$ for all $\omega'_n$, if we do apply $\Psi$ (so $m(\Psi \omega)=-1$) followed by one of these sign-reversing substitutions $\sigma$, we see
\[M(\sigma(\Psi \omega))\geq -m(\Psi \omega)+E(a_1)-1 \geq 1+E(a_1)-1 \geq 0,\] so we may always apply Proposition \ref{proposition - arithmetic one step} without worrying about the possible error of one.
\end{proof}
\end{proposition}

\begin{corollary}[\cite{MR2515388}, Theorem 1, case $k=2$]
\label{corollary - heaviness criterion}
We have $S_n(\theta) \geq 0$ for all $n \geq 0$ if and only if $x(\theta)=0$.
\begin{proof}
By viewing the ergodic sums as an additive cocycle, for all $n>0$ we have $S_n(\theta)=S_{n+1}(0)-1$, so we have by Proposition \ref{proposition - arithmetic on substitutions}:
\[S_{|\Omega'_n|-1}(\theta)=0, \quad M_{|\Omega'_n|-1}(\theta) = \sum_{\substack{i \leq n \\ p_i=0}}E(a_1(\theta_i)), \quad m_{|\Omega'_n|-1}(\theta) = - \sum_{\substack{i \leq n \\ p_i=1}}E(a_1(\theta_i)).\]  So $S_n(\theta)\geq 0$ for all $n$ if and only if $p_i=0 \bmod 2$ for all $i$ such that $\theta_i < 1/2$, which is equivalent to $p_i=0 \bmod 2$ for all $i$.  A direct inductive argument shows that $p_i=0$ for all $i$ if and only if $a_{2i-1}(\theta)=0 \bmod 2$ by considering the action of $g$ \eqref{eqn - g}, which corresponds by Proposition \ref{proposition - heavy theta have zero} to $x(\theta)=0$.
\end{proof}
\end{corollary}
\begin{remark}
Using that $\sigma$ are all homomorphisms, a more constructive version of \eqref{eqn - arbitrary substitution point} is
\[\omega_0 \sigma^{(1)}(\omega_1) \sigma^{(2)}(\omega_2) \cdots \sigma^{(n)}(\omega_n) \cdots,\] which allows a more direct way of computing the word through successive computation of the words $\omega_n$ (given the starting point $x$).
\end{remark}
\begin{lemma}\label{lemma - matrix growth of words}
We always have
\[\left|\sigma^{(n)}(A)\right|=\left|\sigma^{(n)}(B)\right|,\]
and if we define the matrices $M_i=M(\theta_i)$ according to Table \ref{table - matrices for return times}, then
\[M_{n-1} M_{n-2} \cdots M_1 M_0 \left[ \begin{array}{c} 1 \\ 1\end{array}\right]=\left[\begin{array}{c} |\sigma^{(n)}(A)| \\ |\sigma^{(n)}(C)| \end{array}\right].\]
\begin{proof}
The first claim follows directly from the following observation: for all substitutions $\sigma$, the words $\sigma(A)$ and $\sigma(B)$ are always of the same length and always contain the same number of letters drawn from $\{A,B\}$.  That is, within
\[\left(\varphi_{n-1}^{-1} \circ \cdots \circ \varphi_0^{-1}\right)(A \cup B) \subset \tilde{I}_n\] the return time under $R_{\theta_0}$ to $\tilde{I}_n$ is constant, and similarly on the pullback of $C$.  One need only count the number of $C$ and $\{A,B\}$ within $\sigma_n(C)$ and $\sigma_n(\{A,B\})$ to construct the relevant matrices.
\end{proof}
\end{lemma}

\begin{table}[bht]
\begin{tabular}{| c | c |}
\hline Case & $M(\theta)$\\
\hline \hline $a_1(\theta)=0 \bmod 2$, $a_3(\theta) \neq 1$ & $\left[ \begin{array}{c c}(a_1-1)a_2+1 &a_2 \\ (a_1-1)a_2+a_1&a_2+1 \end{array}\right]$\\
\hline $a_1(\theta)=0 \bmod 2$, $a_3(\theta) = 1$ & $\left[ \begin{array}{c c} (a_1-1)a_2+a_1 &a_2+1 \\(a_1-1)a_2+1 &a_2 \end{array}\right]$\\
\hline $a_1(\theta)=1 \bmod 2$, $\neq 1$ & $\left[ \begin{array}{c c} a_1-1& 1 \\ 1& 0 \end{array}\right]$\\
\hline $a_1(\theta)=1$ & $\left[ \begin{array}{c c} 1& 0 \\ 0 &1\end{array}\right]$\\
\hline
\end{tabular}
\caption{The matrices $M(\theta)$ used to determine return times in the induced systems.}\label{table - matrices for return times}
\end{table}

\begin{lemma}\label{lemma - bound on lengths of Omega'(n)}
\[|\Omega_n| \leq |\Omega'_n| \leq |\Omega_{n+1}|.\]
\begin{proof}
The lower inequality is direct in light of \eqref{eqn - sigma' pull off one letter}, recalling that $(\omega'_n)_1=A$.  The upper bound follows from Lemma \ref{lemma - matrix growth of words}, noting that while $\omega'_n$ may or may not be a left factor of $\sigma_n(A)$, it does contain the same number of $\{A,B\}$ versus $C$ as a proper left factor of $\sigma_n(A)$.  Furthermore, the only substitutions for which $|\sigma(C)|>|\sigma(A)|$ are those corresponding to $a_1=0\bmod 2$, $a_3 \neq 0$; such substitutions are \textit{not} followed by $\Psi$.  That is,
\[ \left| \sigma'^{(n)}(A)\right| \leq |\Omega_n|,\] completing the proof of the upper bound.
\end{proof}
\end{lemma}

\begin{figure}\centering
\subfigure[$N=5$, $\sigma(A)=AACAC$]{
\includegraphics[width=2 in]{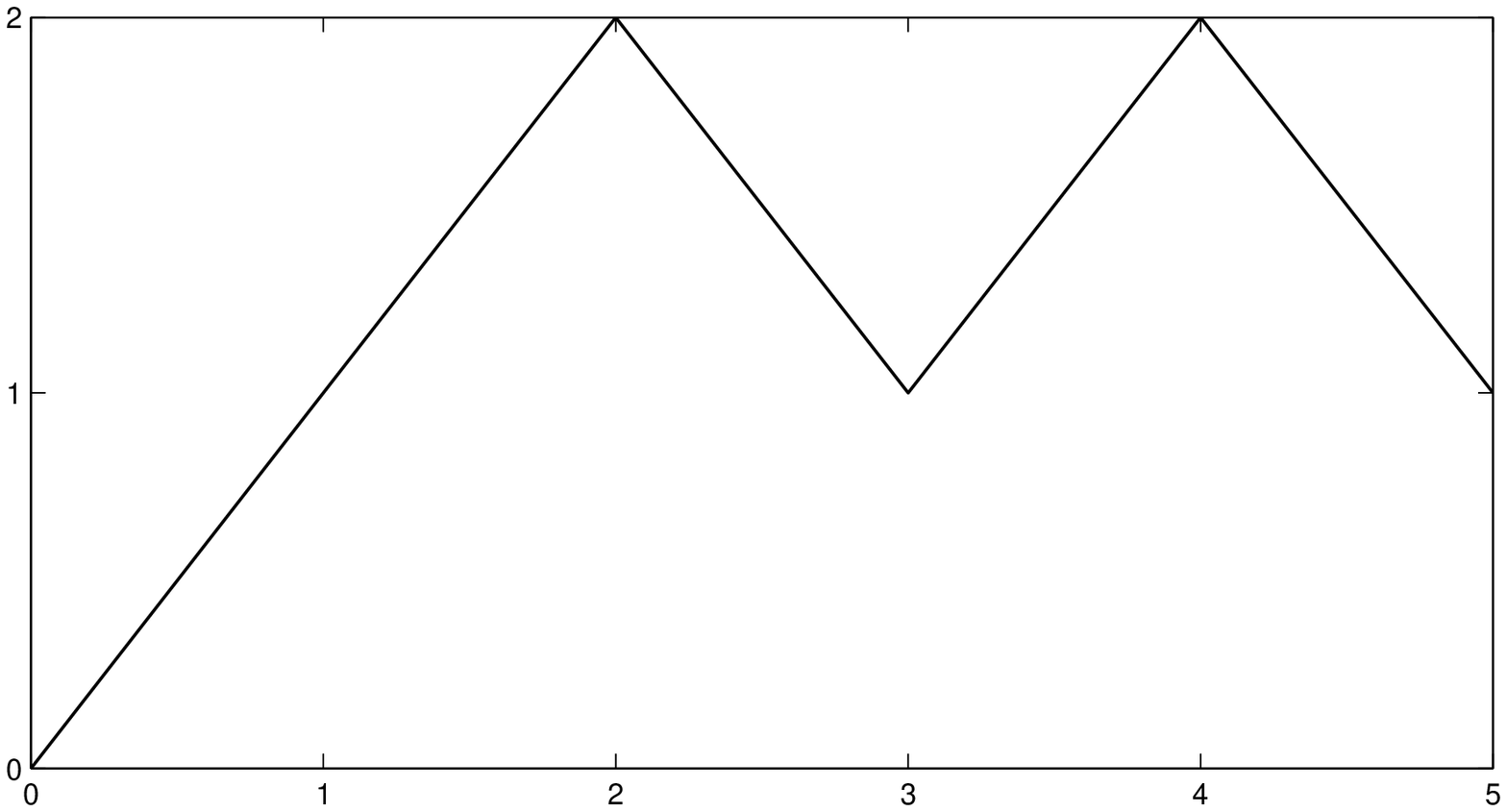}
}
\subfigure[$N=29$, $\sigma^2(A)$]{
\includegraphics[width=2 in]{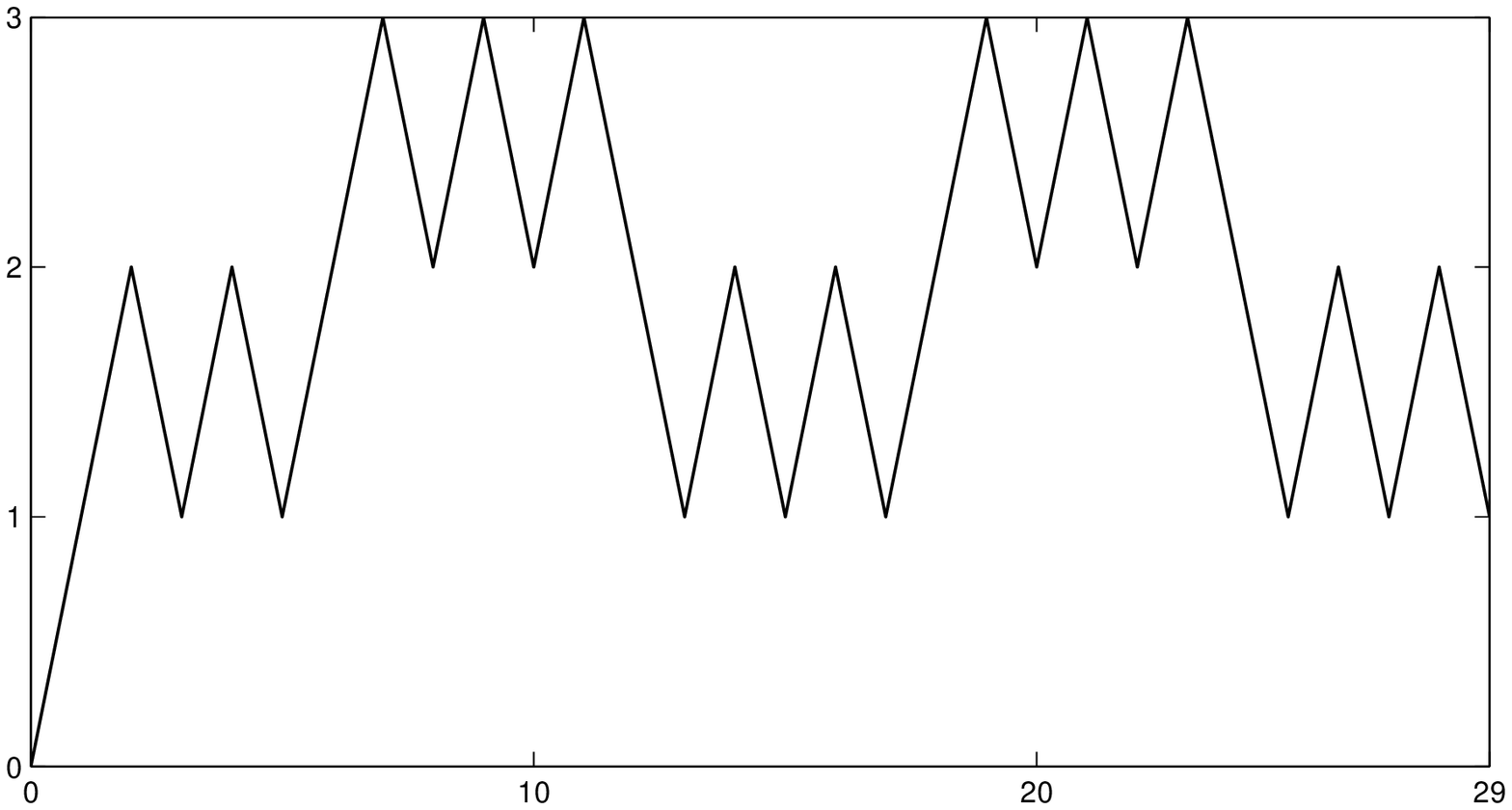}
}
\subfigure[$N=169$, $\sigma^3(A)$]{
\includegraphics[width=2 in]{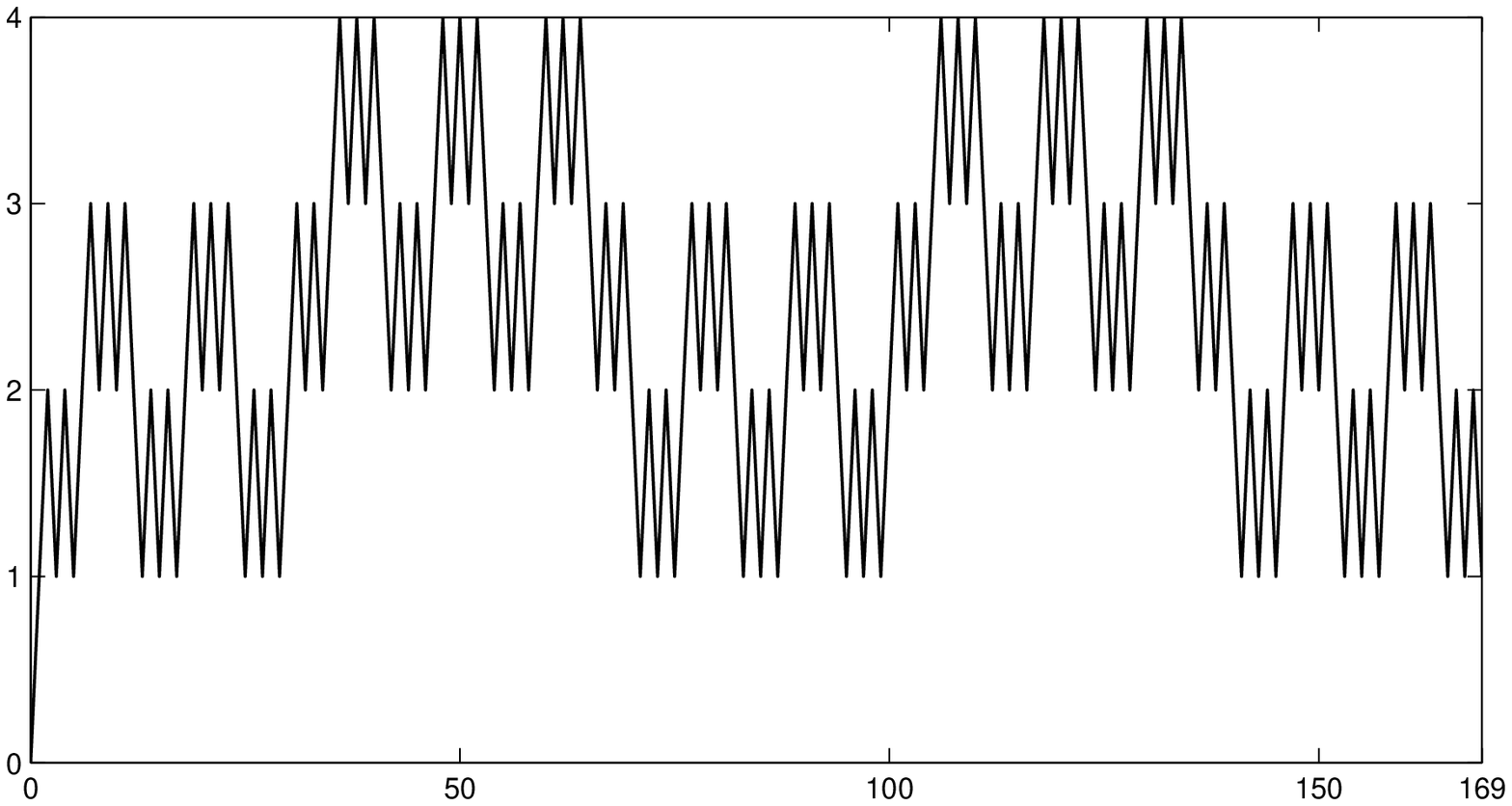}
}
\subfigure[$N=33461$, $\sigma^6(A)$]{
\includegraphics[width=2 in]{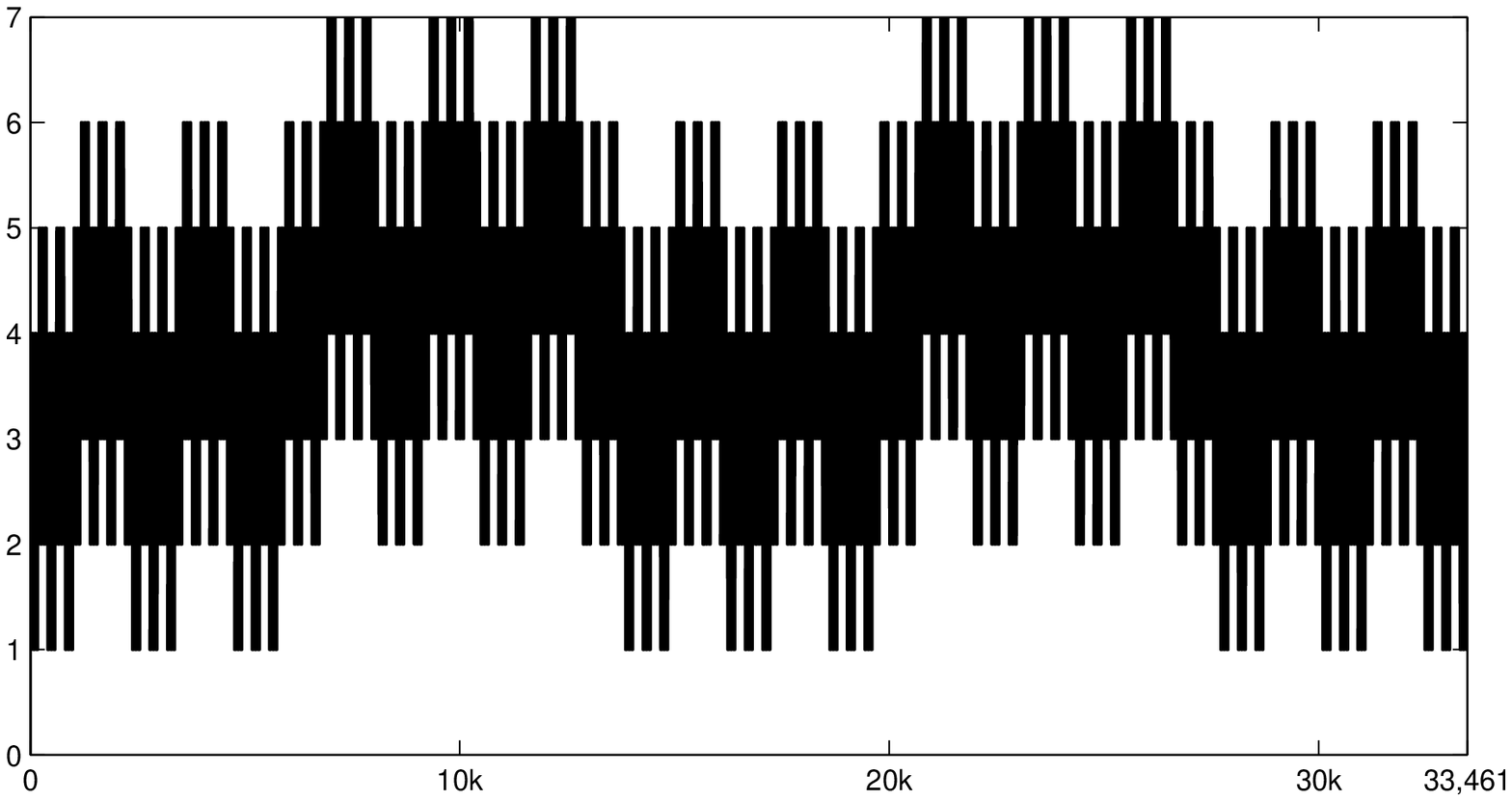}
}
\caption{Plots of $S_i(0)$ for different ranges of $0\leq i \leq N$, where $\theta=\sqrt{2}-1$.}\label{figure - root two example}
\end{figure}

\begin{example}
\label{example - from hensley}
Let $\theta = \sqrt{2} \mod 1 = [2,2,2,\ldots]$.  Then as $\theta$ is a quadratic irrational, the sequence of substitutions $\sigma_i$ is eventually periodic by Proposition \ref{proposition - quadratics are periodic}.  As $g(\theta)=\theta$, the sequence of substitutions is periodic with period one, given by
\[\sigma:  \left\{\begin{array}{l} A \rightarrow AACAC \\ B \rightarrow ABCAC \\ C \rightarrow ABCACAC \end{array} \right.\]
The point $x(\theta)=0$ by Proposition \ref{proposition - heavy theta have zero}, so applying Theorem \ref{theorem - substitution sequence}, the orbit of zero is given by the sequence
\[\lim_{n \rightarrow \infty} \sigma^n(A) = AACACAACACABCACACAACACABCACAC \ldots\]
The self-similar structure of the sequence of ergodic sums $S_n(0)$ is not exact (as $\sigma(B) \neq \sigma(C)$), but nonetheless highly regular.  This regularity was noticed by D. Hensley in \cite[Figure 3.4]{MR2351741}.  We give several plots of $S_n(0)$ for different values of $n$ in Figure \ref{figure - root two example}.  This same self-similarity for developing the orbit of $x(\theta)$ will be seen for any quadratic irrational $\theta$ in light of Proposition \ref{proposition - quadratics are periodic}.
\end{example}

\begin{figure}\centering
\subfigure[$x=0$, with orbit $CACABCACAAC\ldots$]{
\includegraphics[width=5 in]{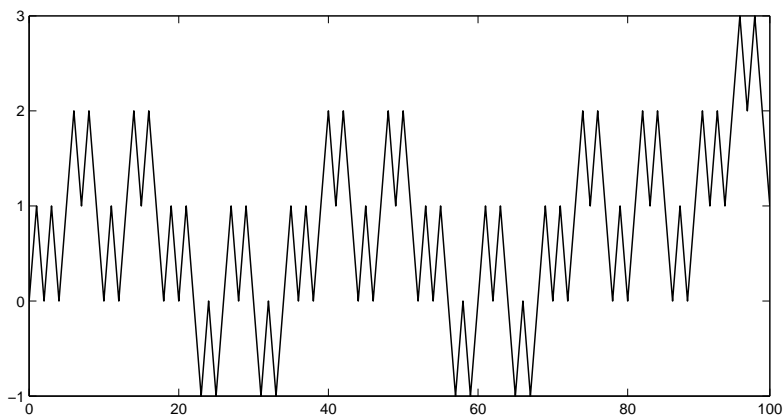}
}
\subfigure[$x=x(\theta)=1/(2\theta)$, with orbit $ABCACAACAC\ldots$]{
\includegraphics[width=5 in]{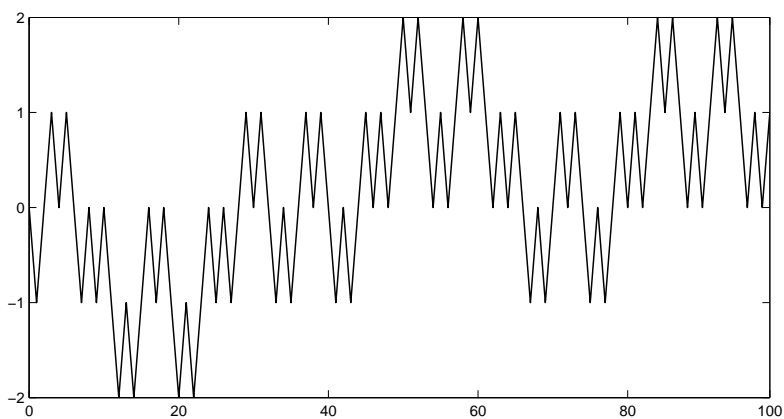}
}
\caption{Plots of $S_i(x)$ for $0 \leq i \leq 100$, where $\theta$ is the golden mean for the two given values of $x$.  Note that as $\theta>1/2$, we have $A\rightarrow -1$, $B,C \rightarrow +1$.}\label{figure - golden mean pictures}
\end{figure}

For quadratic irrational $\theta \notin H$, computation of the point $x(\theta)$ is not too difficult:

\begin{example}\label{example - golden mean}
Let $\theta=[1,1,\ldots]$ be the golden mean.  Recall that $S^1$ will be partitioned such that $A=[(1/2)^+,1^-]$ as $\theta>1/2$.  As $g^2(\theta)=\theta$, and $a_1=1$ corresponds to the identity substitution, the only non-identity substitution generated is
\[\sigma:\left\{ \begin{array}{l} A \rightarrow ABCAC \\ B \rightarrow AACAC \\ C \rightarrow AAC \end{array} \right.\]
So, the orbit of $x(\theta)$ is given by
\[\lim_{n \rightarrow \infty}\sigma^n(A) = ABCACAACACAACABCACAAC\ldots,\]
while the orbit of $0$ is given by
\[\Psi ( \sigma ( \ldots \Psi (AAC)))= CACABCACAACABCACAACAC\ldots.\]

To compute the point $x(\theta)$, we need to determine the intervals $\tilde{I_n}$.  For those $\theta_n=[2,1,1,\ldots]$ we have
\[\delta_n = 1-2\theta_n = 1-2(1-\theta)=2\theta-1.\]  Denote this quantity by $\delta$ for convenience.  For this particular $\theta$ we do not ever have two consecutive $\theta_n<1/2$, so the intervals $I'_{n+1} \subset I_n$ strictly alternate between $[0^+,\delta^-]$ and $[(1-\delta)^+,1^-]$ (for those $n=0\bmod 2$; for odd $n$ we have $\theta_n>1/2$ and $I'_{n+1}=I_n$).  So the sequence of preimages $\tilde{I}_n$ (recall again \eqref{eqn - tilde I}) is given by
\[ \left[0^+, 1^-\right], \quad \left[(1-\delta)^+ , 1^-\right], \quad \left[(1-\delta)^+ ,(1-\delta+\delta^2)^- \right], \ldots\]
whose intersection is given by the geometric series
\[x(\theta) = \sum_{i=0}^{\infty} (-1)^i\delta^i = \frac{1}{1+(2\theta-1)}=\frac{1}{2\theta}.\]
See Figure \ref{figure - golden mean pictures} for both of these orbits.
\end{example}

One particularly striking corollary of Proposition \ref{proposition - arithmetic on substitutions} is the following, which does not seem to be apparent from any other technique:
\begin{corollary}
If $\theta$ is a quadratic irrational, then
\[\lim_{n \rightarrow \infty} \frac{M_n(0)}{|m_n(0)|} \in \mathbb{Q}^*,\] where $\mathbb{Q}^*=\mathbb{Q}\cup\{\infty\}$, and $p/0=\infty$ for any positive integer $p$.  If $\theta_n=\theta_{n+k}$ is a minimal period under the orbit of $g$ and $p_{n+k}=p_{n}+1$, then the ratio tends to one. Furthermore, for any nonnegative $p/q \in \mathbb{Q}^*$, there is a quadratic irrational $\theta$ such that the above ratio has limit $p/q$.
\begin{proof}
We have already shown that $g^n(\theta)$ is eventually periodic for such $\theta$ in Proposition \ref{proposition - quadratics are periodic}.  It follows from Proposition \ref{proposition - arithmetic on substitutions} that $M_n(0)$ and $m_n(0)$ see a periodic sequence of adjustments by bounded integer amounts, which must therefore have rational limit.  If one period reflects a change in the parity of $p$, it will always be followed by the mirrored changes in $M_n$, $m_n$, producing a limit of one.

To produce quadratic irrationals with the desired limit, if $q=0$ then $\theta \in H$ will suffice ($m_n(0) \equiv 1$, and $M_n(0)$ must therefore diverge), and for $p=0$ any $\theta$ such that $a_1(\theta)=1$ and $g(\theta) \in H$ will suffice (here $M_n(0) \equiv 1$).  For $p/q$ with neither zero, just set
\[\theta = [2p,1,1,2q-1,1,1,2p-1,1,1,2q-1,1,1,\ldots],\] and verify that we will first add $p$ to $M_n(0)$, then subtract $q$ from $m_n(0)$, etc.
\end{proof}
\end{corollary}

\section{Proof of Theorem \ref{theorem - aribtrary growth}}

Let $c_n$ and $d_n$ be divergent monotone sequences in $o(n)$ with bounded differences $\Delta c_n$, $\Delta d_n$; we will construct a dense set of $\theta$ such that
\[\limsup_{n \rightarrow \infty} \frac{M_n(0)}{c_n}=\limsup_{n \rightarrow \infty} \frac{|m_n(0)|}{d_n}=1.\]

Any irrational $\theta$ is completely determined by its sequence of partial quotients, which is equivalent to its orbit under $g$, and its orbit under $g$ is completely determined by the sequence of values
\[a_1(\theta_i) \quad (a_1 = 1 \bmod 2), \qquad a_1(\theta_i), \, a_2(\theta_i) \quad (a_1 = 0\bmod 2).\]
Suppose, then, that the first finitely many partial quotients of $\theta$ are prescribed, such that the first $n$ values of $\theta_i$ are fixed.  Without loss of generality, insert an additional single term if necessary so that $p_n=0$ (recall \eqref{eqn - p_n}).  We are now completely free to choose $k$ to construct $\omega'_n$ (refer to \eqref{eqn - omega' words for zero orbit}).  If we denote
\[M(\Omega'_n)=M, \quad m(\Omega'_n)=m, \quad |\Omega'_n|=L_n,\]
it follows from Proposition \ref{proposition - arithmetic on substitutions} that once we choose $k$, we will have
\[M(\Omega'_{n+1})=M+k, \quad m(\Omega'_{n+1})=m.\]  Denote by $L_{n+1}(k)=|\Omega'_{n+1}|$ as a function of $k$.

Assume first that $M < c_{L_n}$, so we wish to increase the maximal sum compared to the sequence $c_n$.  Then let $a_1(\theta_n)$ be odd, so
\[\omega'(n+1) = A^{k+1} B^{k}.\]
From \eqref{eqn - sigma' pull off one letter} and the previous observation that $|\sigma^{(n)}(A)|=|\sigma^{(n)}(B)|$, it follows that
\[L_{n+1}(k) = |\tilde{\omega}| + 2k |\sigma^{(n+1)}(A)|,\]
where
\[\tilde{\omega} = \sigma'^{(n+1)}(A).\]

Consider, then, the proper left factors $A^{i}$ of $\omega'(n+1)$ for $i=1,2,\ldots, k+1$.  Applying Proposition \ref{proposition - arithmetic one step}, the new maximal sum $M+k$ is achieved at a time $N$, where
\[|\tilde{\omega}|+(k-1)|\sigma^{(n)}(A)| \leq N \leq |\tilde{\omega}|+k|\sigma^{(n)}(A)|.\]  As $c_n \in o(n)$, we may choose $k \geq 1$ to be minimal such that
\[\frac{M+k}{c(|\tilde{\omega}|+k|\sigma^{(n)}(A)|)} \geq 1.\]

If, however, we had $M \geq c_{L_n}$, then we would wish to not greatly increase $M$ compared to $c_n$.  In this case, let $\theta_n=[2,k,1,\ldots]$, and pass directly to considering the word
\[\sigma'^{(n+1)}(C)=\sigma'^{(n)}(A^{k+1}B^{k-1}C),\] as $C$ is always a left factor of $\omega'_{n+1}=\Psi(\omega'_{n+2})$ in this case.  Then the maximal sum reached for this word is $M+1$, but its length is (similarly to before)
\[L_{n+1}(k) = |\tilde{\omega}| + 2k |\sigma^{(n)}(A)|.\]  We are now in the position of being able to increase the length of the word \textit{without} increasing the maximal sum of $M+1$, so as $c_n$ is divergent, choose $k \geq 1 $ minimal such that
\[\frac{M+1}{c(|\tilde{\omega}|+k|\sigma^{(n)}(A)|)} \leq 1.\]

After applying $g$ twice (to skip past the next $\theta_k>1/2$), then, we find ourselves able to manipulate the growth of the minimal sums $m(n)$.  Continuing in this fashion, then, we construct a dense set of $\theta$ (as the initial string of partial quotients was arbitrary).  That the $\limsup$s are actually one follows from the minimal choice of $k$ and that $\Delta c_n$, $\Delta d_n$ are bounded.

To prove the analogous statements where one of $M_n$, $m_n$ is desired to remain bounded, one need only repeat the same arguments using $\theta_n \in H$ (recall \eqref{eqn - heavy theta}) so that the value $p_n$ is eventually constant. \qedhere

The statement of Theorem \ref{theorem - aribtrary growth} applies as well to $M_n(x(\theta))$ and $m_n(x(\theta))$; the proof is simpler, in fact, as the map $\Psi$ is not a concern, and the possible error of one from Proposition \ref{proposition - arithmetic on substitutions} is not an asymptotic concern.  This process is highly amenable to diagonalization techniques.  For example:
\begin{corollary}
Given a countable collection of sequences $c^{(i)}_n$ and $d^{(i)}_n$, all of which are divergent and in $o(n)$, such that
\[c^{(1)}_n \leq c^{(2)}_n \leq \ldots, \quad d^{(1)}_n \geq d^{(2)}_n \geq \ldots,\] there is a dense set of $\theta$ for which
\[c^{(i)}_n \in o(M_n(0)), \quad |m_n(0)| \in o(d^{(i)}_n)\] for all $i$.
\begin{proof}
Apply Theorem \ref{theorem - aribtrary growth} after using a diagonalization process to construct $c_n$, $d_n$, both monotone, divergent, and in $o(n)$ such
\[c^{(i)}_n \in o(c_n), \quad d_n \in o(d^{(i)}_n).\qedhere\]
\end{proof}
\end{corollary}

\begin{figure}\centering
\subfigure[$\theta$ exhibiting very slow growth of $M_n(0)$; this portion of the graph will repeat $2^{16}$ times with no additional growth.]{
\includegraphics[width=5 in]{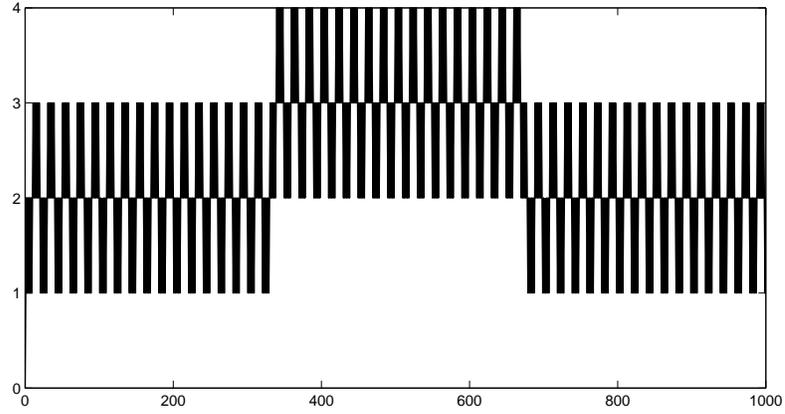}
}
\subfigure[$\gamma(\theta)$ exhibiting very fast growth of $M_n(0)$; this sawtooth pattern will continue to climb by repeating itself $E(2^{16})/2$ times.]{
\includegraphics[width=5 in]{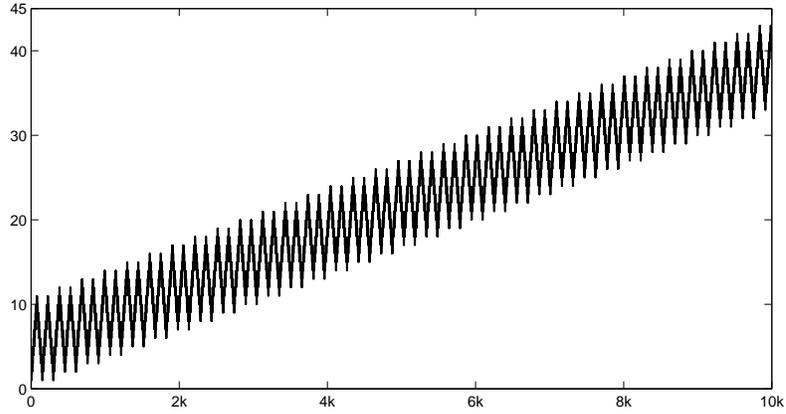}
}
\caption{Two different extreme growth rates for $\theta$ and $\gamma(\theta)$.}\label{figure - two extreme growths}
\end{figure}

Many permutations of the above corollary are possible.  For example, we may construct a dense set of $\theta$ such that the discrepancy sums grow in both directions faster than any $n^{1-\epsilon}$ (but necessarily in $o(n)$, of course!), or such that the discrepancy sums are bounded below, but $M_n(0)$ grows slower than all iterated logarithms (but necessarily divergent, of course!), etc.  See Figure \ref{figure - two extreme growths} for an example where for both $\theta$ and $\gamma(\theta)$ we have $m_n \geq 1$, but $M_n(\theta)\notin o(n^{1-\epsilon})$ for any $\epsilon>0$ while $M_n(\gamma(\theta)) \in o (\log^{(i)}n)$ for all $i$.  In Figure \ref{figure - two extreme growths} we set
\[\theta = [2,2^2,2,2^{2^2},2,2^{2^{2^2}},2,\ldots].\]

Using diagonalization techniques one may similarly find a dense set of $\theta$ such that
\[\limsup_{i \rightarrow \infty} \frac{M_{n_i(j)}(0)}{c^{(j)}_{n_i(j)}}=1\] for an arbitrary collection of divergent sequences $c^{(j)}_n$ in $o(n)$ for different subsequences $n_i(j) \rightarrow \infty$ depending on $j$, and similarly for the $|m_n(0)|$ and a collection of sequences $d^{(j)}_n$.

Truly, beyond the constraints of \eqref{eqn - unbounded but not linear}, any asymptotic behavior desired is possible.

\section{Proof of Theorem \ref{theorem - arbitrary interval for ratio}}
Suppose that
\begin{equation}\label{eqn - for arbirary ration theorem}\liminf_{n \rightarrow \infty} \frac{M_n(0)}{|m_n{0}|}=r_1, \quad \limsup_{n \rightarrow \infty} \frac{M_n(0)}{|m_n(0)|}=r_2.\end{equation}

That the set of accumulation points of the sequence is the entire closed interval $[r_1,r_2]$ is direct and is left to the reader.  Let an arbitrary finite string of partial quotients $a_1,\ldots,a_N$ be given which determine $\theta_i$ for $i=0,1,\ldots,n-1$, and for convenience again assume without loss of generality that $p_n=0$.

Now let $c_n$ and $d_n$ be arbitrary integer-valued strictly increasing sequences such that $\Delta c_n$ and $\Delta d_n$ are in $O(1)$ and
\[\liminf_{n \rightarrow \infty} \frac{c_n}{d_n}= \rho_1, \quad \limsup_{n \rightarrow \infty} \frac{c_n}{d_n}= \rho_2.\]  Furthermore, assume that
$c_1 > M(\Omega'_n)=M$ and $d_1 > |m(\Omega'_n)|=m$.

Continue the continued fraction expansion of $\theta$ in the following way:
\[\theta_n=[2(c_1-M)+1, 2(d_1-m),2(c_2-c_1),2(d_2-d_1),\ldots].\]
Then $\Omega'_n$ will see the sequence of $M(\Omega'_{n+2k})=c_k$ and $m(\Omega'_{n+2k})=-d_{k}$; the bounded differences $\Delta c_n$ and $\Delta d_n$ ensure that the limiting behavior is the same as the limiting behavior along the subsequence of times $|\Omega'_n|$. \qedhere

\begin{figure}\centering
\includegraphics[width=5 in]{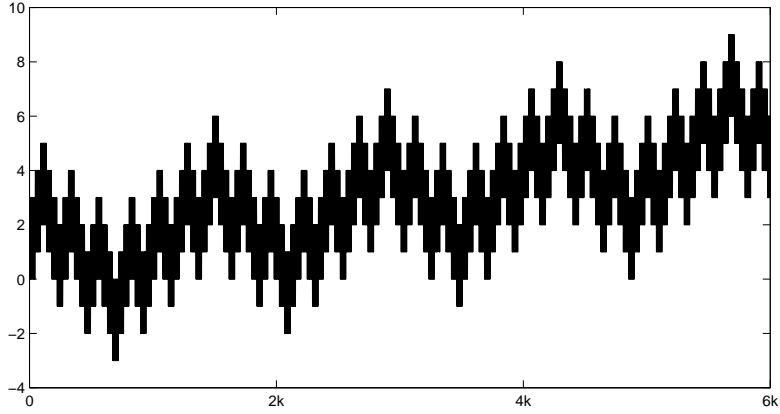}
\caption{A specific $\theta$ for which $M_n(0)/|m_n(0)|$ has limit two; refer to \eqref{eqn - prdicted strangegrowth} and note the changes to $M$, $m$.}
\label{figure - ratio two}
\end{figure}

\begin{example}\label{example - ratio two}
Suppose that $\theta=[1,2,3,4,\ldots]$.  Then we begin computing the sequence of values $M_n(0)$ and $|m_n(0)|$ according to Proposition \ref{proposition - arithmetic on substitutions}:
\begin{equation}\label{eqn - prdicted strangegrowth}\begin{array}{|c |c |c |c|}
\hline  \theta_0 = [1,2,3,4,\ldots] & p=0  &E(a_1)=0 & (M,|m|)=(1,1)\\
 \theta_1 = [3,3,4,5,\ldots] & p=1 & E(a_1)=1 & (M,|m|)=(1,0)\\
 \theta_2 = [1,3,4,5,\ldots] & p=1 & E(a_1)=0 & (M, |m|)=(1,0)\\
 \theta_3 = [4,4,5,6,\ldots] & p=0 & E(a_1)=2 & (M,|m|)=(3,0)\\
 \theta_4 = [5,6,7,8,\ldots] & p=0 & E(a_1)=2 & (M,|m|)=(5,0)\\
 \hline \theta_5 = [1,6,7,8\ldots] & p=0 & E(a_1)=0 & (M,|m|)=(5,0)\\
 \theta_6 = [7,7,8,9,\ldots] & p=1 & E(a_1)=3 & (M, |m|)=(5,3)\\
\vdots & \vdots & \vdots & \vdots\\
\hline
\end{array}\end{equation}
The pattern is seen to continue in groups of five terms.  Over the terms $\theta_{5k}$ through $\theta_{5k+4}$, we will subtract $2k+1$ from $m$ while adding $2(2k+2)$ to $M$.  We therefore have $\rho_1=\rho_2=2$, or
\[\lim_{n \rightarrow \infty} \frac{M_n(0)}{|m_n(0)|}=2.\]
See Figure \ref{figure - ratio two} for this $\theta$.
\end{example}

\section{Proof of Theorem \ref{theorem - bounded type are log}}

\begin{lemma}\label{lemma - pass to all points for bounded type}
Suppose that $f(x)$ is a step function on $S^1$ with $k < \infty$ discontinuities, and denote $V(f)$ the variation of $f$.  Define $S_n(x)$, $M_n(x)$ and $m_n(x)$ as before.  As we have not restricted $f$ to be integer-valued, define
\[\rho_N(x) = \left(M_N-m_N\right)(x).\]
Let $n$ be such that $q_n \leq N < q_{n+1}$. Then for any $x,y \in S^1$:
\[\rho_N(y) \leq \rho_{q_{n+2}}(x)+a_{n+1}V(f).\]
\begin{proof}
Consider the set $\{x+i \theta\}$ for $i=0,1,\ldots,q_{n}-1$.  Choose $0 \leq j < q_{n}$ such that $x+j \theta$ is closest to $y$.  Then the distance between $x+j\theta$ and $y$ is no larger than $q_n^{-1}$.  For each discontinuity $d_i$ there are therefore at most $a_{n+1}$ preimages of $d_i$ within this interval for time $L=0,1,\ldots,q_{n+1}-1$.  It follows that $f(x+(j+i)\theta)=f(y+i\theta)$ for all but at most $k\cdot a_{n+1}$ of $i=0,1,\ldots,N < q_{n+1}$.  As $j+i$ is less than $q_n+q_{n+1} \leq q_{n+2}$, the lemma follows.
\end{proof}
\end{lemma}

Assume that $a_i(\theta) \leq M$ for all $i$.  Then (continuing with existing notation) we see that for some $C>1$ independent of $\theta$
\begin{equation}\label{eqn - bounded type time growth}C^{\frac{n-1}{2}} \leq |\Omega'_n| \leq (M+1)^{2n+2}.\end{equation}  The lower bound is due to the exponential decay in the length of the interval $\tilde{I}_n$ (any $C<2$ eventually suffices, as $\tilde{I}_{n+1}$ is less than half as large as $\tilde{I}_n$ at least half the time, with the $n-1$ accounting for the possibility that $I'_1=I_0$, or $\theta_0>1/2$).  The upper bound follows from Lemma \ref{lemma - matrix growth of words}, Lemma \ref{lemma - bound on lengths of Omega'(n)}, and the bound $a_i(\theta)\leq M$.
while at the same time,
\begin{equation}\label{eqn - bounded type sum growth}\frac{n-1}{2} \leq \rho_{|\Omega'_n|}(0) \leq \frac{nM}{2};\end{equation}
the lower inequality is due to the fact that at most half of the words $\Omega'_n=\Omega'_{n+1}$ (corresponding to those $\theta_n>1/2$) and for the rest, $\rho(\Omega_{n+1})\geq \rho(\Omega_n)+1$, as $E(a_1)\geq 1$ for these $\theta_n<1/2$.  The upper bound follows as $E(a_i(\theta))\leq M/2$ for all $i$.

Now, for any $N$ let $k$ be chosen such that
\[|\Omega_k| \leq N \leq |\Omega_{k+1}|.\]
From \eqref{eqn - bounded type time growth}:
\[k C_1\leq \log |\Omega'_k| \leq \log(N) \leq \log |\Omega'_{k+1}| \leq kC_2,\] for two constants $C_1$ and $C_2$ which do not depend on $k$.  From \eqref{eqn - bounded type sum growth}:
\[\frac{(k+1)M}{2} \geq \rho_{|\Omega'_{k+1}|}(0) \geq \rho_N(0) \geq \rho_{|\Omega'_k|}(0) \geq \frac{k-1}{2},\] so $\rho_n(0) \sim \log (n)$.  The full theorem now follows from Lemma \ref{lemma - pass to all points for bounded type}.
\section*{Acknowledgements}
The author is greatly indebted to many people for support and helpful conversations over the development of this paper.  The original impetus for studying this problem came from a problem posed by M. Boshernitzan while the author was a Ph.D. student at Rice University, while a rudimentary form of Theorem \ref{theorem - aribtrary growth} arose from discussions at PRIMA 2008 during a visit supported by the University of New South Wales.  The author is currently supported by the Center for Advanced Studies at Ben Gurion University of the Negev, where Barak Weiss has provided invaluable suggestions on improving the clarity of an early draft.  Of course, any mistakes or unclear passages in the current form are entirely the author's responsibility.

\bibliographystyle{plain}
\bibliography{bibfile}
\end{document}